[1994/12/01]
\documentclass[draft]{article}
\title{Orthogonal maximal abelian *-subalgebras of the $6\times 6$ matrices }
\author{Kyle Beauchamp\\ and\\Remus Nicoara\\  Vanderbilt University, Nashville TN}

\date{}

\usepackage{amsmath,amsthm}

\chardef\bslash=`\\ 

\usepackage{amsfonts}
\newtheorem{thm}{Theorem}[section]

\newtheorem{lem}[thm]{Lemma}
\newtheorem{obs}[thm]{Remark}

\theoremstyle{definition}

\theoremstyle{remark}

\newcommand{\eval}[2][\right]{\relax
  \ifx#1\right\relax \left.\fi#2#1\rvert}

\let\abs=\envert

\newcommand{\HermDiagThreeC}{\begin{pmatrix}
1 & 1 & 1 & 1 & 1 & 1 \cr
1 & -1 & \bar{a} & -\bar{a}  &\bar{b}  & -\bar{b}   \cr
1 & a & -1 & \bar{c} & -a  & -\bar{c}  \cr
1 & -a  & c & -1  & \bar{\alpha} &\bar{\beta}  \cr
1 & b  & -\bar{a} & \alpha & 1& \bar{\gamma}  \cr
1 & -b &-c  & \beta  & \gamma & 1 \cr
\end{pmatrix}}

\newcommand{\HermDiagThreeA}{\begin{pmatrix}
1 & 1 & 1 & 1 & 1 & 1 \cr
1 & -1 & \bar{a} & \bar{b}  &-\bar{a}  & -\bar{b}   \cr
1 & a & -1 & -a & \bar{c}  & -\bar{c}  \cr
1 & b  & -\bar{a} & -1  & \bar{\alpha} &\bar{\beta}  \cr
1 & -a  & c & \alpha & 1& \bar{\gamma}  \cr
1 & -b &-c  & \beta  & \gamma & 1 \cr
\end{pmatrix}
}

\newcommand{\HermDiagThreeB}{\begin{pmatrix}
1 & 1 & 1 & 1 & 1 & 1 \cr
1 & -1 & \bar{a} & \bar{b}  &-\bar{a}  & -\bar{b}   \cr
1 & a & -1 & \bar{c} & -a  & -\bar{c}  \cr
1 & b  & c & -1  & \bar{\alpha} &\bar{\beta}  \cr
1 & -a  & -\bar{a} & \alpha & 1& \bar{\gamma}  \cr
1 & -b &-c  & \beta  & \gamma & 1 \cr
\end{pmatrix}
}
\newcommand{\HermDiagThreeD}{\begin{pmatrix}
1 &  1  & 1        &         1 &        1     & 1             \cr
1 & -1  & \bar{a}  & -\bar{a}  & \bar{b}      & -\bar{b}      \cr
1 &  a  & -1       & -a        &      \bar{c} & -\bar{c}      \cr
1 & -a  & -\bar{a} & -1        & \bar{\alpha} &\bar{\beta}    \cr
1 &  b  & c        & \alpha    & 1            & \bar{\gamma}  \cr
1 & -b  & -c       & \beta     & \gamma       & 1             \cr
\end{pmatrix}
}

\newcommand{\HermDiagThreeE}{\begin{pmatrix}
1 &  1   & 1        &         1 &        1     & 1             \cr
1 & -1   & \bar{a}  & \bar{b}  & -\bar{a}      & -\bar{b}      \cr
1 &  a   & -1       & \bar{c}  &   -\bar{c}    & -a      \cr
1 & b    & c        & -1        & \bar{\alpha} &\bar{\beta}    \cr
1 &  -a  & -c       & \alpha    & 1            & \bar{\gamma}  \cr
1 & -b   & -\bar{a} & \beta     & \gamma       & 1             \cr
\end{pmatrix}
}

\newcommand{\HermDiagThreeEb}{\begin{pmatrix}
1 &  1           & 1              &         1 &        1      & 1             \cr
1 & -1           & \bar{x}        & -y        & -\bar{x}      & y      \cr
1 &  x           & -1             & t         &         -t    & -x      \cr
1 & -\bar{y}     & \bar{t}        & -1        & \bar{\alpha}  &\bar{\beta}    \cr
1 &  -x          & -\bar{t}       & \alpha    & 1             & \bar z  \cr
1 & \bar{y}      & -\bar{x}       & \beta     & z       & 1             \cr
\end{pmatrix}
}

\newcommand{\HermDiagFive}{\begin{pmatrix}
1 &  1           & 1              &         1 &        1      & 1             \cr
1 & -1           & \bar{a}        & \bar{b}   & -\bar{a}      & -\bar{b}      \cr
1 &  a           & -1             & \bar{x}   &    \bar{y}    & \bar{z}      \cr
1 & b            & x              & -1        & \bar{\alpha}  &\bar{\beta}    \cr
1 &  -a          & y              & \alpha    & -1             & \bar{\gamma}  \cr
1 & -b           & z              & \beta     & \gamma        & -1             \cr
\end{pmatrix}
}

\newcommand{\HermDiagFour}{\begin{pmatrix}
1 &  1           & 1              &         1 &        1      & 1             \cr
1 & -1           & \bar{a}        & \bar{b}   & -\bar{a}      & -\bar{b}      \cr
1 &  a           & -1             & \bar{x}   &    \bar{y}    & \bar{z}      \cr
1 & b            & x              & -1        & \bar{\alpha}  &\bar{\beta}    \cr
1 &  -a          & y              & \alpha    & -1             & \bar{\gamma}  \cr
1 & -b           & z              & \beta     & \gamma        & 1             \cr
\end{pmatrix}
}

\newcommand{\HermDiagTwo}{\begin{pmatrix}
1 &  1           & 1              &         1 &        1      & 1             \cr
1 & -1           & \bar{a}        & \bar{b}   & -\bar{a}      & -\bar{b}      \cr
1 &  a           & -1             & \bar{x}   &    \bar{y}    & \bar{z}      \cr
1 & b            & x              & 1        & \bar{\alpha}  &\bar{\beta}    \cr
1 &  -a          & y              & \alpha    & 1             & \bar{\gamma}  \cr
1 & -b           & z              & \beta     & \gamma        & 1             \cr
\end{pmatrix}
}

\newcommand{\HermDiagOne}{\begin{pmatrix}
1 &  1           & 1              &         1 &        1      & 1             \cr
1 & -1           & \bar{a}        & \bar{b}   & -\bar{a}      & -\bar{b}      \cr
1 &  a           & 1             & \bar{x}   &    \bar{y}    & \bar{z}      \cr
1 & b            & x              & 1        & \bar{\alpha}  &\bar{\beta}    \cr
1 &  -a          & y              & \alpha    & 1             & \bar{\gamma}  \cr
1 & -b           & z              & \beta     & \gamma        & 1             \cr
\end{pmatrix}}

\newcommand{\CircNine}
{\ensuremath{
\begin{pmatrix}
   
   1 & 1 & 1 & 1 & 1 & 1 & 1 & 1 & 1 \cr 1 & 
    -1 & -1 & y & y^{-1} & -y^{-3} & y^
   {-3} & y & y^{-1} \cr 1 & y & y^
   3 & -y^2 & y^3 & y & -1 & 
    -y & y \cr 1 & y & y & -1 & y^{-1} & y^{-3} & -y^{-3} & -1 & y^{-1} \cr 1 & -y & -y^2 & y^3 & y^3 & 
    -1 & y & y & y \cr 1 & y & y^3 & y^
   3 & y & y^{-1} & y^{-1} & y & y^
   2 \cr 1 & y^3 & y & y & y^2 & y^{-1} & y^{-1} & y^3 & y \cr 1 & 
    -y^2 & -y & y & y & -1 & y & y^3 & y^
   3 \cr 1 & y^3 & y & -y & y & y & -1 & 
    -y^2 & y^3 \cr  
        \end{pmatrix}
}}
\newcommand{\nineten}{
\ensuremath{
\begin{pmatrix}
1 & 1 & 1 & 1 & 1 & 1 & 1 & 1 & 1 \cr 1 &
   -1 & \epsilon^3 & \epsilon^3 & -1 & \epsilon^9 & \epsilon^8 & \epsilon^7 & \epsilon \cr 1& \epsilon^4 & -1 & \epsilon^7 & \epsilon & \epsilon^3 & -1& \epsilon^9 & \epsilon^9 \cr 1 & \epsilon^3 & \epsilon^7 & -1 & \epsilon & \epsilon^8 & \epsilon^9 & \epsilon^3 & -1 \cr 1 &
 \epsilon^9 & \epsilon & -1 & -1 & \epsilon^3 & \epsilon^7 & \epsilon^2 & \epsilon^7 \cr
 1 & \epsilon^9 & -1 & \epsilon & \epsilon^3 & -1 & \epsilon & \epsilon^7 & \epsilon^6 \cr 1 & \epsilon & \epsilon^{7} & \epsilon^9 & \epsilon^ 6 & \epsilon & -1 & -1 & \epsilon^3 \cr 1 & \epsilon^7 & \epsilon^9 & \epsilon^4 & \epsilon^9 & -1 & \epsilon^3 & -1 & \epsilon \cr 1 & -1 & \epsilon^2 & \epsilon^9 & \epsilon^7 & \epsilon^7 & \epsilon^3 & \epsilon & -1 \cr 
    \end{pmatrix}
}}

\newcommand{\romeo}{
$\left(
\begin{array}{ccccccccccc}
1 & 1 & 1 & 1 & 1 & 1 & 1  & 1  & 1  & 1  & 1 \cr
1 & y & y & a & b & b & -1 & x  & -1 & x  & -1 \cr
1 & y & b & b & a & y & x  & x  & -1 & -1 & -1 \cr
1 & a & b & y & y & b & -1 & -1 & x  & -1 & x \cr
1 & b & a & y & b & y & x  & -1 & -1 & -1 & x \cr
1 & b & y & b & y & a & -1 & -1 & x  & x  & -1 \cr
1 & -1 & x & -1 & x & -1 & - \frac{x}{a}   & -\frac{x}{b}   & - \frac{x}{y}   & - \frac{x}{y}  & - \frac{x}{b}  \cr
1 & x & x & -1 & -1 & -1 & -\frac{x}{b}   & - \frac{x}{y}   & - \frac{x}{b} & - \frac{x}{y}   & -\frac{x}{a}    \cr
1 & -1 & -1 & x & -1 & x & - \frac{x}{y} & - \frac{x}{b}  & - \frac{x}{b}   & -\frac{x}{a}   & - \frac{x}{y}   \cr
1 & x & -1 & -1 & -1 & x & - \frac{x}{y} & - \frac{x}{y}   & -\frac{x}{a}   & - \frac{x}{b}   & - \frac{x}{b}  \cr
1 & -1 & -1 & x & x & -1 & -\frac{x}{b}   & - \frac{x}{a}   & - \frac{x}{y}& - \frac{x}{b}   & -\frac{x}{y}    \cr 
\end{array}
\right)$
}

\newcommand{\liketau}{
\ensuremath{
\begin{pmatrix}
  1 & 1 & 1 & 1 & 1 & 1 & 1 & 1 & 1 & 1 \cr 
   1 & 1 & a^{-2} & a^{-1} & a^{-2} & a^
   2 & 1 & a^2 & 1 & a \cr 1 & a^2 & 1 & a^
   {-2} & a^{-2} & 1 & a & a^2 & 1 & a^{-1} \cr 1 & 1 & 1 & 1 & a^{-2} & a^2 & a^
   2 & a & a^{-2} & a^{-1} \cr 1 & a^
   2 & a^{-2} & 1 & 1 & a & 1 & a^2 & a^
   {-2} & a^{-1}\cr 1 & 1 & a^{-1} & a^{-2} & 1 & 1 & a^2 & a^2 & a^
   {-2} & a \cr 1 & a^2 & 1 & a^{-2} & a^{-1} & a^2 & 1 & 1 & a^{-2} & a \cr 1 & a^
   2 & a^{-2} & 1 & a^{-2} & 1 & a^
   2 & 1 & a^{-1}& a \cr 1 & a & a^
   {-2} & a^{-2} & 1 & a^2 & a^
   2 & 1 & 1 & a^{-1} \cr 1 & a &a^{-1} &a^{-1} & a^{-1}& a & a & a & a^{-1}& 1 \cr 
   \end{pmatrix}
   }
   }

\begin{document}
\maketitle

\renewcommand{\sectionmark}[1]{}
\begin{abstract}
We construct new pairs of orthogonal maximal abelian $*$-subalgebras of $M_6(\mathbb C)$, by classifying all self-adjoint complex Hadamard matrices of order 6. In particular, we exhibit a non-affine one-parameter family of non-equivalent Hadamard matrices of order 6. In the last part of the paper we present other previously unknown examples of complex Hadamard matrices of higher orders.
\end{abstract}
\section{Introduction}$ $

Let $(M_n(\mathbb C),Tr)$ denote the algebra of $n\times n$ complex matrices with the usual trace. For $X\in M_n(\mathbb C)$, denote by $X^*$ the conjugate transpose of $X$. A subalgebra $\mathcal A$ of $M_n(\mathbb C)$ is called a MASA if it is maximal abelian and closed under the $*$ operation. It is easy to see that $\mathcal A$ is a MASA if and only if it is unitarily conjugate to the algebra $\mathcal D_n$ of $n\times n$ diagonal matrices, i.e. $\mathcal A=U\mathcal D_nU^*$ for some $U\in M_n(\mathbb C)$ unitary matrix. 

We say that two MASA's $\mathcal A_1,\mathcal A_2$ are \emph{orthogonal} if $\mathcal A_1\cap \mathcal A_2=\mathbb C$ and the vector subspaces $\mathcal A_1\ominus\mathbb C$ and $\mathcal A_2\ominus\mathbb C$ are orthogonal, with respect to the inner product $<x,y>=Tr(y^*x)$, $x,y\in M_n(\mathbb C)$. This is equivalent to saying that the square of inclusions:
$$\mathfrak{C}=\left(\begin{matrix}
\mathcal A_1 & \subset{} & M_n(\mathbb{C}) \cr
\cup   &           & \cup \cr
\mathbb{C} & \subset{} & \mathcal A_2
\end{matrix},Tr\right)$$
is a commuting square in the sense of \cite{Po1},\cite{Po2} (see also \cite{GHJ}). 

We may assume, up to unitary conjugacy of commuting squares, that $\mathcal A_1=\mathcal D_n$ and $\mathcal A_2=U\mathcal D_nU^*$, for some unitary $U\in M_n(\mathbb C)$. In this notation, the orthogonality of $\mathcal A_1,\mathcal A_2$ amounts to $U$ having all entries of the same absolute value $\frac{1}{\sqrt n}$, hence $H=\sqrt{n} U$ is a \emph{complex Hadamard matrix}. Thus, being given a commuting square $\mathfrak{C}$ of this form is equivalent to having a complex $n\times n$ Hadamard matrix.

Recall that two complex Hadamard matrices are \emph{equivalent} if there exist unitary diagonal matrices $D_1,D_2$ and permutation matrices $P_1,P_2$ such that
$$H_2=P_1D_1H_1D_2P_2$$

It is easy to see that equivalence of Hadamard matrices corresponds to isomorphism of commuting squares, via the identification described above. Our interest in Hadamard matrices, mainly in obtaining one-parameter families of non-equivalent Hadamard matrices, comes from the possibility of constructing subfactors from the corresponding commuting squares (see for instance \cite{Jo2}). However, it is hard to decide if such subfactors are non-isomorphic, or to compute their principal graphs.

Besides their connections to von Neumann algebras (\cite{Haagerup},\cite{HJ},\cite{Jo2},\cite{MW},\cite{Ni},\\ \cite{Petrescu}, \cite{Po2}), complex Hadamard matrices have numerous other applications such as the theory of error correcting codes (\cite{CH}), spectral sets and Fuglede's conjecture (\cite{T}). They play a very important role in quantum information theory, in the construction of teleportation and dense coding schemes (\cite{We}). 

While all Hadamard matrices of orders up to 5 are classified (\cite{Haagerup}), it seems very hard to describe Hadamard matrices of higher orders, such a classification not being known even for $n=6$. For $n$ composite, some constructions of parametric families of Hadamard matrices whose entries are linear functions (also called \emph{affine} families, see \cite{TZ}) are presented in \cite{Di},\cite{MRS}. There is however no general procedure of constructing such families with non-affine entries, or for $n$ prime.  A catalogue of most known complex Hadamard matrices of small order (up to order 16) can be found in \cite{TZ}. 

In this paper we classify, up to equivalence, all Hadamard matrices $H$ of order 6 that are self-adjoint, i.e. $H=H^*$, where $H^*$ denotes the conjugate transpose of $H$. We thus obtain a new one-parameter non-affine family:
$$H(\theta)=\begin{pmatrix}
1 & 1 & 1 & 1 & 1 & 1 \cr
1 & -1 & \bar{x} &-y  &-\bar{x}  &y   \cr
1 & x & -1 & t & -t  & -x  \cr
1 &-\bar{y}  & \bar{t} & -1  & \bar{y} &-\bar{t}   \cr
1 &-x  & -\bar{t} & y & 1 &\bar{z}  \cr
1 & \bar{y} &-\bar{x}  &-t  & z & 1 \cr
\end{pmatrix}$$
where $\theta\in [-\pi,-arcos(\frac{-1+\sqrt{3}}{2})]\cup [arcos(\frac{-1+\sqrt{3}}{2}),\pi]$ and the variables $x,y,z,t$ are given by: 
$$y=exp(i\theta),\text{ } z=\frac{1+2y-y^2}{y(-1+2y+y^2)}$$
$$x=\frac{1+2y+y^2-\sqrt{2}\sqrt{1+2y+2y^3+y^4}}{1+2y-y^2}$$
$$t=\frac{1+2y+y^2-\sqrt{2}\sqrt{1+2y+2y^3+y^4}}{-1+2y+y^2}$$

Our family is non-affine, in particular it is not obtained by modifying linearly the entries of a tensor product of $2\times 2$ and $3\times 3$ Hadamard matrices. Such constructions result in subfactors having intermediate subfactors, thus of non-trivial first relative commutant. We analysed computationally our family $H(\theta)$ and we conjecture that it yields subfactors of principal graph $A_\infty$.

In \cite{Ni} the second author introduced a condition for commuting squares, called \emph{the span condition}, which is sufficient to imply isolation of a commuting square in the class of commuting squares (up to isomorphisms). In particular, applying this to commuting squares arising from Hadamard matrices, we obtained a sufficient condition for isolation of a Hadamard matrix among all Hadamard matrices (up to equivalence). In \cite{TZ} the notion of \emph{defect} of a Hadamard matrix was introduced. Saying that the defect of a matrix is zero is equivalent to the span condition.

It is not settled whether the span condition is also necessary for isolation. In $\cite{TZ}$ a possible counter-example is provided: is is shown that no affine family stems from the Bjorck-Froberg 'cyclic 6 roots' matrix $C_6^{(0)}$ (\cite{Bj},\cite{Haagerup}), while its defect is non-zero. Therefore, it is asked whether this matrix is isolated among all Hadamard matrices. We answer negatively to this question, by showing that the family $H(\theta)$ contains a matrix equivalent to $C_6^{(0)}$. Thus, the question of whether isolation is equivalent to the span condition remains open.

In the last part of the paper we present other new examples of complex Hadamard matrices of dimensions $9,10,11$. These examples were found using computers, by doing a numerical search for Hadamard matrices satisfying certain symmetry conditions.

We would like to thank Teodor Banica, Ingemar Bengtsson, Dietmar Bisch, Wes Camp, Romeo Maciuca, Wojciech Tadej and Karol Zyczkowski for fruitful discussions and correspondence. Kyle Beauchamp was supported in part by NSF under Grant No. DMS 0353640 (REU Grant), and Remus Nicoara was supported in part by NSF under Grant No. DMS 0500933.

\section{Self-adjoint Hadamard matrices of order 6}

In this section we classify, up to equivalence, all complex self-adjoint Hadamard matrices of order 6. We prove the following theorem, stating that there exists a non-affine one-parameter family of such matrices.

\begin{thm}\label{main} Let $H\in M_6(\mathbb C)$ be a self-adjoint Hadamard matrix. Then $H$ is equivalent to $H(\theta)$, for some $\theta\in [-\pi,-arcos(\frac{-1+\sqrt{3}}{2})]\cup [arcos(\frac{-1+\sqrt{3}}{2}),\pi]$, where:

$$H(\theta)=\begin{pmatrix}
1 & 1 & 1 & 1 & 1 & 1 \cr
1 & -1 & \bar{x} &-y  &-\bar{x}  &y   \cr
1 & x & -1 & t & -t  & -x  \cr
1 &-\bar{y}  & \bar{t} & -1  & \bar{y} &-\bar{t}   \cr
1 &-x  & -\bar{t} & y & 1 &\bar{z}  \cr
1 & \bar{y} &-\bar{x}  &-t  & z & 1 \cr
\end{pmatrix}$$
and the parameters $x,y,z,t$ are given by: 

$$y=exp(i\theta),\text{ } z=\frac{1+2y-y^2}{y(-1+2y+y^2)}$$
$$x=\frac{1+2y+y^2-\sqrt{2}\sqrt{1+2y+2y^3+y^4}}{1+2y-y^2}$$
$$t=\frac{1+2y+y^2-\sqrt{2}\sqrt{1+2y+2y^3+y^4}}{-1+2y+y^2}$$

\end{thm}

\begin{obs} In \cite{TZ} it is asked whether the Bjorck-Froberg 'cyclic 6 roots' matrix:
$$C_6^{(0)}=
\begin{pmatrix}
1 & 1 & 1 & 1 & 1 & 1 \cr
1 & -1 & -d & -d^2 & d^2 & d  \cr
1 & -\bar{d} & 1 & d^2 & -d^3 & d^2  \cr
1 & -\bar{d^2} & \bar{d^2} & -1  & d^2 & -d^2  \cr
1 & \bar{d^2} & -\bar{d^3} & \bar{d^2} & 1 & -d \cr
1 & \bar{d} & \bar{d^2} & -\bar{d^2} & -\bar{d} & -1 \cr
\end{pmatrix}, d=\frac{1-\sqrt{3}}{2}+i{}(\frac{\sqrt{3}}{2})^\frac{1}{2}$$
is isolated among complex Hadamard matrices. It is known that no affine Hadamard family stems from $C_6^{(0)}$ (\cite{TZ}). However, this matrix does not satisfy the span condition we introduced in \cite{Ni}, or equivalently its defect is non-zero, in the sense of \cite{TZ}. If $C_6^{(0)}$ were isolated, it would follow that the span condition, which is sufficient to ensure isolation, is not necessary for isolation. However, Theorem \ref{main} shows that there exists a continuum of non-equivalent Hadamard matrices containing $C_6^{(0)}$, since $H(\theta_0)$ is equivalent to $C_6^{(0)}$ for $\theta_0=2Arg(d)$. Indeed: $$P H(\theta_0) P^{-1}= C_6^{(0)}$$
 where $P$ is the permutation matrix:
$$P=\begin{pmatrix}
   1 & 0 & 0 & 0 & 0 & 0 \cr 0 & 1 & 0 & 0 & 0 & 0 \cr 0 & 0 & 0 & 0 & 1 
& 0 \cr 0 & 0 & 0 & 1 & 0 & 0 \cr 0 &
   0 & 0 & 0 & 0 & 1 \cr 0 & 0 & 1 & 0 & 0 & 0 \cr 
\end{pmatrix}$$
In particular, the question of wether isolation is equivalent to the span condition remains open.

\end{obs}
\noindent Recall that a complex Hadamard matrix $H=(h_{k,l})\in M_n(\mathbb C)$ is said to be \emph{dephased} or in \emph{normal form} if $h_{1,k}=h_{k,1}=1$ for all $k=1,..,n$. The next lemma shows that, in order to classify all self-adjoint Hadamard matrices, one only needs to look at dephased self-adjoint Hadamard matrices.

\begin{lem}
Let $H\in M_n(\mathbb C)$ be a self-adjoint Hadamard matrix. Then $H$ is equivalent to a dephased self-adjoint Hadamard matrix.
\end{lem}
\begin{proof}
Since $H=(h_{k,l})_{1\leq k,l\leq n}$ is hermitian, $h_{k,k}$ are real and thus they belong to $\{-1,1\}$. We may assume, by eventually multiplying $H$ by $-1$, that $h_{1,1}=1$. Consider the matrix $H'=(h_{k,l}\bar h_{k,1}\bar h_{1,l})_{1\leq k,l\leq n}$. $H'$ is clearly equivalent to $H$ and $ h'_{k,1}= h'_{1,l}=1$ for all $1\leq k,l \leq n$. Morover, $H'$ is hermitian since $\bar h'_{k,l}=\bar h_{k,l}h_{k,1}h_{1,l}=h_{l,k}\bar h_{1,k}\bar h_{l,1}=h'_{l,k}$.

\end{proof}

\noindent We now recall two easy lemmas involving algebraic manipulations of complex numbers.

\begin{lem}\label{3vectors}

If $x, y, z \in \mathbb{C}$ such that $\abs{x}=\abs{y}=\abs{z}=1$ and $x+y+z=0$, then $x=z\epsilon$ and $y=z\epsilon^2$, where $\epsilon\in
\{\frac{1}{2}+\frac{i\sqrt{3}}{2},\frac{1}{2}-\frac{i\sqrt{3}}{2}\}$.
\end{lem}
\begin{proof}
By conjugating $x+y+z=0$ we obtain $\frac{1}{x}+\frac{1}{y}+\frac{1}{z}=0$. Solving and eliminating $x$ yields $\frac{1}{y+z}=\frac{1}{y}+\frac{1}{z}$. Equivalently, $(\frac{y}{z})^2+\frac{y}{z}+1=0$, which shows that $y=z\epsilon$, with $\epsilon$ as above.  Since $1+\epsilon+\epsilon^2=0$, $x=-y-z=-z(1+\epsilon)=z\epsilon^2$.
\end{proof}
\begin{lem}\label{4vectors} If $x,y,z,t \in \mathbb{C}$ such that $\abs{x}=\abs{y}=\abs{z}=\abs{t}=1$ and $x+y+z+t=0$, then $x\in\{-y,-z,-t\}$.
\end{lem}

\begin{proof}
We have: $(x+y)(x+z)(x+t)=x^2(x+y+z+t)+xyz+xyt+xzt+yzt=xyz+xyt+xzt+yzt=xyzt (\bar x+\bar y+\bar z+ \bar t)=0.$
\end{proof}

\noindent The following lemma is also used in \cite{Haagerup}, towards the classification of complex Hadamard matrices of order 5 .
 
\begin{lem}\label{HLemma} Let $u,v,s,t$ be complex numbers on the unit circle. Then: $$(u+v)(\bar s+\bar t)(\bar u s+\bar v t)\in \mathbb R$$
\end{lem}

\begin{proof}
$(u+v)(\bar s+\bar t)(\bar u s+\bar v t)=(u\bar s+v\bar t+u\bar t+v\bar s)(\bar u s+\bar v t)=2+(u\bar v \bar s t+\bar u v s \bar t)+(\bar u v+u\bar v)+(\bar s t+s\bar t)$ is real, since $z+\bar z$ is real for every $z\in\mathbb C$.

\end{proof}

\noindent We now proceed with the proof of Theorem \ref{main}. Since $H$ is hermitian, its diagonal elements belong to $\{-1,1\}$. Morover, since for every permutation matrix $P$ the matrix $H$ is equivalent to $PHP^{-1}$, and $PHP^{-1}$ is still hermitian, it is enough to consider the following six possibilities for the diagonal of $H$:
 
$ $

\noindent $Diag(H)\in$ $\{(1,1,1,1,1,1), (1,-1,1,1,1,1),(1,-1,-1,1,1,1),(1,-1,-1,-1,1,1)$\\ $(1,-1,-1,-1,-1,1), (1,-1,-1,-1,-1,-1)\}$. 

$ $

\noindent We start by showing that the diagonal of $H$ can not be $(1,1,1,1,1,1)$. This is the most difficult of the six cases we need to analyse. Indeed, in all the other cases the existence of a $1$ and a $-1$ on one of the rows of $H$ will allow us to apply  Lemma \ref{4vectors}, thus reducing the number of variables.

\begin{lem}
\textbf{(a).} Let $H$ be a complex $6\times 6$ Hadamard matrix of the form: $$H=\begin{pmatrix} 1 & 1 & 1 & 1 & 1 & 1\\
1 & 1 & x & \bar y & . & . \\
1 & \bar x & 1 & z & . & . \\
1 & y & \bar z & 1 & . & . \\
1 & . & . & . & . & . \\
1 & . & . & . & . & . 

\end{pmatrix}$$
Then two of $x,y,z$ must be equal.

\noindent\textbf{(b).} Let $H\in M_6(\mathbb C)$ be a self-adjoint, dephased, complex Hadamard matrix. Then the diagonal of $H$ can not be $(1,1,1,1,1,1)$. 
\end{lem}
\begin{proof}

(a). Assume, by contradiction, $x\not=y\not=z\not=x$. Denote the last two elements on the second and third rows of $H$ by $u=h_{2,5},v=h_{2,6}$, respectively $s=h_{3,5},t=h_{3,6}$. Using the orthogonality of the first three rows of $H$ we obtain:
$$2+x+\bar y = -(u+v)$$ $$2+x+\bar z =-(\bar s + \bar t)$$ $$ 1+2\bar x +yz=-(\bar u s+\bar v t)$$ Lemma \ref{HLemma} implies:
\begin{equation}\label{eqx}(2+x+\bar y)(2+x+\bar z)(1+2\bar x +yz)\in\mathbb R\end{equation}

\noindent The same argument applied to rows 1,2,4 respectively rows 1,3,4, shows:

\begin{equation}\label{eqy}(2+y+\bar z)(2+y+\bar x)(1+2\bar y +zx)\in\mathbb R\end{equation}
and
\begin{equation}\label{eqz}(2+z+\bar x)(2+z+\bar y)(1+2\bar z +xy)\in\mathbb R\end{equation}

\noindent Expanding the product in (\ref{eqx}) and using $x\bar x=y\bar y=z\bar z=1$, we obtain: 

$ $

\noindent $(x^2+\bar y \bar z + 4+x\bar y + x\bar z+4x+2\bar y +2\bar z)(1+2\bar x+yz)=x^2yz+x^2+(2\bar x\bar y\bar z+4xyz)+x\bar y+x\bar z+(4\bar x \bar y+xy) + (4\bar x  \bar z +xz)+ (\bar y\bar z+4yz) + (8\bar x + 6x)+(4\bar y +2y)+(4\bar z+2z)+13\in \mathbb R$. 

$ $

\noindent Since $(2\bar x\bar y\bar z+2xyz)+(\bar x\bar y +xy)+(\bar x\bar z+xz)+4(\bar y\bar z+yz)+6(\bar x +x)+(2\bar y +2y)+(2\bar z+2z)+13\in\mathbb R$, by substracting it from the previous expression it follows:
$$x^2yz+x^2+2xyz+x\bar y+x\bar z+3\bar x \bar y+3\bar x  \bar z -3\bar y\bar z+2\bar x+2\bar y+2\bar z\in \mathbb R$$
Thus: $$x^2yz+x^2+x(\bar x+\bar y+\bar z)-6\bar y \bar z + (2xyz+3\bar x \bar y+3\bar x  \bar z +3\bar y\bar z+2\bar x+2\bar y+2\bar z)\in\mathbb R$$
Let $$S=2xyz+3\bar x \bar y+3\bar x  \bar z +3\bar y\bar z+2\bar x+2\bar y+2\bar z$$
Using that $(6yz+6\bar y\bar z)$ is real yields:
$$x^2yz+x^2+x(\bar x+\bar y+\bar z)+6yz + S\in \mathbb R$$
Similarly, by expanding (\ref{eqy}) and reducing real terms we obtain: 
$$y^2zx+y^2+y(\bar x+\bar y+\bar z)+6zx + S\in\mathbb R$$
The number S is the same, since (\ref{eqy}) is just a circular permutation $(x,y,z)\rightarrow (y,z,x)$ of (\ref{eqx}), and the formula for $S$ is invariant to permutations of $x,y,z$. Substracting the two previous expressions and cancelling S, we obtain: 
$$xyz(x-y)+x^2-y^2 +(x-y)(\bar x+\bar y+\bar z)-6(x-y)z\in \mathbb R$$
Thus: $$(x-y)(xyz+x+y+\bar x+\bar y+\bar z-6z)\in \mathbb R$$
 Hence:
$$(x-y)(xyz+x+y+z+\frac{1}{x}+\frac{1}{y}+\frac{1}{z}-7z)\in\mathbb R$$
Using that a complex number is real iff it equals its conjugate, we obtain:

$$(x-y)(xyz+x+y+z+\frac{1}{x}+\frac{1}{y}+\frac{1}{z}-7z)=(\frac{1}{x}-\frac{1}{y})(\frac{1}{xyz}+\frac{1}{x}+\frac{1}{y}+\frac{1}{z}+x+y+z-\frac{7}{z})$$
Multiplying the previous equality by $-\frac{xy}{x-y}$ yields: 
$$-xy(xyz+x+y+z+\frac{1}{x}+\frac{1}{y}+\frac{1}{z}-7z)=\frac{1}{xyz}+\frac{1}{x}+\frac{1}{y}+\frac{1}{z}+x+y+z-\frac{7}{z}$$ 
Thus: $$\frac{7}{z}=xy(xyz+x+y+z+\frac{1}{x}+\frac{1}{y}+\frac{1}{z})-7xyz+\frac{1}{xyz}+\frac{1}{x}+\frac{1}{y}+\frac{1}{z}+x+y+z$$
 Let $$T=-7xyz+\frac{1}{xyz}+\frac{1}{x}+\frac{1}{y}+\frac{1}{z}+x+y+z$$ $$R= xyz+x+y+z+\frac{1}{x}+\frac{1}{y}+\frac{1}{z}$$
We showed:
$$\frac{7}{z}=xyR+T$$
By repeating the same argument for the relations (\ref{eqy}),(\ref{eqz}), and using $y\not=z$, we obtain:
$$\frac{7}{x}=yzR+T$$
where $R,T$ are the same as in the previous equation, since their formulas are symmetric in $x,y,z$. Substract the last two equations: $$\frac{7}{z}-\frac{7}{x}=(x-z)yR$$
Multiplying by $\frac{xz}{x-z}$, we obtain: $$7=xyzR$$
 This implies $|R|=7$. However, by the triangle inequality we have: $$|R|=|xyz+x+y+z+\frac{1}{x}+\frac{1}{y}+\frac{1}{z}|\leq 1+1+1+1+1+1+1=7$$
 Since $x,y,z$ were assumed distinct, we can't have equality and thus we have reached a contradiction.

$ $

\noindent (b). Reasoning by contradiction, assume that $H$ satisfies the hypothesis. Denoting $H$ as in part (a) we know that two of $x,y,z$ are equal. We analyse the three possible cases: $x=y$, $y=z$, $x=z$.

\textbf{Case I: $x=y$.}  From (\ref{eqx}) we obtain: 
\begin{equation}\label{xz} (2+x+\bar x)(2+x+\bar z)(1+2\bar x+xz)\in\mathbb R\end{equation}
If $x=y=-1$,  using the orthogonality of columns 1,2 we obtain $h_{6,2}=-h_{5,2}$. Orthogonality of rows 5,6 yields $$2h_{6,5}+h_{6,3} \bar h_{5,3}+h_{6,4}\bar h_{5,4}=0$$  By Lemma \ref{4vectors}, we must have $$h_{6,5}=-h_{6,3}\bar h_{5,3}=-h_{6,4} \bar h_{5,4}$$
In particular: $h_{5,4}\bar h_{5,3}=h_{6,4}\bar h_{6,3}$. Using this together with the orthogonality of columns $3,4$ yields: $$1+z+h_{5,4}\bar h_{5,3}=0$$
 which together with Lemma \ref{3vectors} implies $z\in\{\epsilon,\epsilon{}^2\}$, $\epsilon=exp(2\pi i /3)$. It is immediate to check that this contradicts equation (\ref{eqz}).

\noindent This shows that $x\not=-1$. Thus,  $2+x+\bar x\not=0$ and after dividing by it in equation (\ref{xz}) we obtain: 
$$(2+x+\bar z)(1+2\bar x+xz)=4+4\bar x+2x+2xz+2\bar x\bar z+x^2 z+\bar z \in\mathbb R$$
Thus: $2\bar x+x^2 z+\bar z \in\mathbb R$, i.e. $\frac{2}{x}+x^2z+\frac{1}{z}=2x+\frac{1}{x^2z}+z$. After multiplying by $x^2z$ and simplifying: $$(x^2-1)(xz-1)^2=0$$

\noindent We have $x\not=1$, since if $x=y=1$ the orthogonality of columns 1,2 implies $4+h_{5,2}+h_{6,2}=0$, which is impossible. Thus $xz=1$, so $x=y=\bar z$. It is easy to check that in this case relation (\ref{eqz}) holds true. 

\noindent We substract the sum of the elements of column 4 (which is 0) from the sum of the elements of column 3 (which is also 0):
$$-\bar h_{5,4}-\bar h_{6,4}+h_{5,3}+h_{6,3}=0$$
Last equation together with Lemma \ref{4vectors} yields one of three possibilities:

\textbf{I.(i).} $h_{6,4}=-h_{5,4}$ and $h_{6,3}=-h_{5,3}$. In this case, the orthogonality of columns 1 and 3 implies $x=-1$, which we showed it is not possible.

\textbf{I.(ii).}  $h_{6,4}=\bar h_{5,3}$ and $h_{5,4}=\bar h_{6,3}$. From the orthogonality of columns 3,4 we obtain: $$-\frac{(1+x)^2}{2}=h_{5,3}h_{6,3}$$
The last equality implies $x\in\{i,-i\}$. However, the sum of the elements of the third column of $H$ is 0: $$2+2x+h_{5,3}+h_{6,3}=0$$
which contradics the triangle inequality: $2\sqrt{2}=|2+2x|=|h_{5,3}+h_{6,3}|\leq 2$.

\textbf{I.(iii).} $h_{6,4}=\bar h_{6,3}$ and $h_{5,4}=\bar h_{5,3}$. Substracting the inner product of columns 4,2 from the inner product of columns 2,3 we obtain:
$$h_{5,2}\bar h_{5,3}-\bar h_{5,2}\bar h_{5,3}+h_{6,2}\bar h_{6,3}-\bar h_{6,2}\bar h_{6,3}=0$$Applying Lemma \ref{4vectors}, we have three possibilities: 

\textbf{I.(iii).1.}  $h_{5,2}\bar h_{5,3}=\bar h_{5,2}h_{5,3}$ and $h_{6,2}\bar h_{6,3}=\bar h_{6,2}h_{6,3}$. Thus $h_{5,2}\bar h_{5,3}$ and $h_{6,2}\bar h_{6,3}$ are real. Combining this with the orthogonality of columns 2,3 we obtain $x$ real, thus $x=\pm 1$, a contradiction.

\textbf{I.(iii).2.}  $h_{5,2}\bar h_{5,3}=-h_{6,2}\bar h_{6,3}$. Using this relation together with the orthogonality relation between columns 2,3, we obtain $x=-1$, contradiction.

\textbf{I.(iii).3.} $h_{5,2}\bar h_{5,3}=\bar h_{6,2}h_{6,3}$. This equality together with the orthogonality of columns 2,3 implies $\bar x\in \mathbb R$, thus $x=\pm 1$, contradiction.

\noindent This ends the analysis of the case when $x=y$.

\textbf{Case II: $y=z$.} This case can be treated similarly to Case I.

\textbf{Case III: $x=z$.} As in the first case, one of the following holds: $x=1$, $x=-1$, or $x=\bar{y}=z$. However, since the sum of the elements on the third row of $H$ is $0$ we can't have $x=1$. Also, if $x=-1$ then equation (\ref{eqy}) implies $y=\pm 1$, which contradicts the orthogonality of rows 2,4. Thus we must have $$x=\bar{y}=z$$
Writing that the sum of the entries of column 2 (which is 0) equals the sum of the conjugates of the entries of column 4 (also equal to 0), we obtain: $$h_{5,4}+h_{6,4}=\bar h_{5,2}+\bar h_{6,2}$$ Lemma \ref{4vectors} divides now the problem in three cases.

\textbf{Case III.(i).} $h_{6,4}=-h_{5,4}$ and $h_{6,2}=-h_{5,2}$. In this case the orthogonality of rows 1,2 forces $x=-1$, which we showed it is not possible.

\textbf{Case III.(ii).} $h_{5,4}=\bar h_{6,2}$ and $h_{6,4}=\bar h_{5,2}$. This case can be treated similarly to Case I.(ii), but by looking at columns 2,4 instead of columns 3,4.

\textbf{Case III.(iii).}$h_{5,4}=\bar h_{5,2}$ and $h_{6,4}=\bar h_{6,2}$. Again, this case can be treated similarly to Case I.(iii), by substracting the inner product of columns 4,3 from the inner product of the columns 2,3.

\end{proof}

\begin{lem}
Let $H\in M_6(\mathbb C)$ be a self-adjoint, dephased, complex Hadamard matrix. Then the diagonal of $H$ can not be $(1,-1,1,1,1,1)$.
\end{lem}
\begin{proof}
Since the first two columns of $H$ are orthogonal, the sum of the elements on the second column of $H$ is 0. Two of these elements being $-1,1$, the sum of the other four equals 0. By applying Lemma \ref{3vectors}, we may assume, after eventually permuting some rows and the corresponding columns of $H$, that the second column of $H$ is $(1,-1,a,b,-a,-b)$. Thus:
$$H=\HermDiagOne$$

\noindent  Using the ortoghonality of columns 3,5 of $H$, we obtain: $$2y+x\alpha+z\bar \gamma=0$$ Lemma \ref{4vectors} implies $y=-\alpha x$, so  $z= -x\alpha\gamma$. Similarly, considering the ortoghonality of columns 4,6 of $H$, we obtain $\beta=-\alpha\gamma$. After making these substitutions, the orthogonality of columns 3,4 yields $$1+\bar{a}b=0$$ 
while from the orthogonality of columns 4,5 it follows $$1-\bar{a} b=0$$
Last two relations are clearly contradictory.
\end{proof}

\begin{lem}
Let $H\in M_6(\mathbb C)$ be a self-adjoint, dephased, complex Hadamard matrix. Then the diagonal of $H$ can not be $(1,-1,-1,1,1,1)$.
\end{lem}

\begin{proof}
As in the previous lemma, we may assume:
$$H=\HermDiagTwo$$
The orthogonality of columns $4,6$, together with Lemma \ref{4vectors}, implies 
$$\beta=-\alpha\gamma, z=\alpha\gamma x$$
Since the inner product of columns 3,5 is zero, we have: 
$$ z\bar\gamma=\alpha x=-y$$
The orthogonality of columns 3,4 implies: $$2x=1+\bar a b$$
and thus $\bar a b=x=1$, while the orthogonality of columns 3,6 yields: $$2x\alpha\gamma=1-\bar a b=0$$
which is a contradiction.
\end{proof}
\begin{lem}Let $H\in M_6(\mathbb C)$ be a self-adjoint, dephased, complex Hadamard matrix. If the diagonal of $H$ is $(1,-1,-1,-1,1,1)$, then $H$ is equivalent to $H(\theta)$, for some $\theta\in [-\pi,-arcos(\frac{-1+\sqrt{3}}{2})]\cup [arcos(\frac{-1+\sqrt{3}}{2}),\pi]$, where:

$$H(\theta)=\begin{pmatrix}
1 & 1 & 1 & 1 & 1 & 1 \cr
1 & -1 & \frac{1}{x} &-y  &-\frac{1}{x}  &y   \cr
1 & x & -1 & t & -t  & -x  \cr
1 &-\frac{1}{y}  & \frac{1}{t} & -1  & \frac{1}{y} &-\frac{1}{t}   \cr
1 &-x  & -\frac{1}{t} & y & 1 &\frac{1}{z}  \cr
1 & \frac{1}{y} &-\frac{1}{x}  &-t  & z & 1 \cr
\end{pmatrix}$$
and the parameters $x,y,z,t$ are given by: 

$$y=exp(i\theta),\text{ } z=\frac{1+2y-y^2}{y(-1+2y+y^2)}$$
$$x=\frac{1+2y+y^2-\sqrt{2}\sqrt{1+2y+2y^3+y^4}}{1+2y-y^2}$$
$$t=\frac{1+2y+y^2-\sqrt{2}\sqrt{1+2y+2y^3+y^4}}{-1+2y+y^2}$$

\end{lem}
\begin{proof}
By applying lemma \ref{4vectors} to the elements on the second column of $H$, which sum up to 0, it follows that the second column of $H$ has to be of one of the forms:$$\begin{pmatrix}
1\cr
-1\cr
a\cr
-a\cr
b\cr
-b\cr
\end{pmatrix},
\begin{pmatrix}
1\cr
-1\cr
a\cr
b\cr
-a\cr
-b\cr
\end{pmatrix}, \text{ or }
\begin{pmatrix}
1\cr
-1\cr
a\cr
b\cr
-b\cr
-a\cr
\end{pmatrix}$$
We may discard the third option, since it is equivalent with the second option by permuting rows 5,6 and columns 5,6 of $H$. Indeed, this operation does not change the diagonal of $H$. By applying Lemma \ref{4vectors} to the third column of $H$, which contains $\bar a$ since $H$ is hermitian, we obtain the following possibilites for columns 2,3 of $H$:
$$\begin{pmatrix}
1 & 1\cr
-1 &\bar{a} \cr
a & -1\cr
-a & -\bar{a} \cr
b & c \cr
-b & -c \cr
\end{pmatrix},
\begin{pmatrix}
1 & 1\cr
-1 &\bar{a} \cr
a & -1\cr
-a & c \cr
b & -\bar{a} \cr
-b & -c \cr
\end{pmatrix},
\begin{pmatrix}
1 & 1\cr
-1 &\bar{a} \cr
a & -1\cr
-a & c \cr
b & -c \cr
-b & -\bar{a} \cr
\end{pmatrix},
\begin{pmatrix}
1 & 1\cr
-1 &\bar{a} \cr
a & -1\cr
b & -\bar{a} \cr
-a & c \cr
-b & -c \cr
\end{pmatrix},
\begin{pmatrix}
1 & 1\cr
-1 &\bar{a} \cr
a & -1\cr
b & c \cr
-a & -\bar{a} \cr
-b & -c \cr
\end{pmatrix},
\begin{pmatrix}
1 & 1\cr
-1 &\bar{a} \cr
a & -1\cr
b & c \cr
-a & -c \cr
-b & -\bar{a} \cr
\end{pmatrix}$$
We may remove the third arrangement, since it is equivalent to the second by permuting rows 5,6 and columns 5,6 of $H$, and replacing $b$ by $-b$. We now analyse the five cases left.

$ $

\noindent \textbf{Case I.}  $$H=\HermDiagThreeD$$
Using that the inner product of the last two columns of $H$ is $0$ we obtain:
$$-1+\bar\alpha\beta+2\gamma=0$$
From lemma \ref{4vectors} we have:
$$\gamma=1, \alpha=\beta$$
Since the sum of the elements of column 5 of $H$ is $0$, using $\gamma=1$ yields:
$$3+b+c+\alpha=0$$
which implies $b=c=\alpha=-1$. This contradicts the fact that the sum of the elements of column 6 is $0$: $2-b-c+\beta+\gamma=0$.

$ $

\noindent\textbf{Case II.}$$H=\HermDiagThreeC$$
Using the orthogonality of columns 3,4 yields:
$$-2c-\bar a\bar \alpha-c\bar\beta=0$$ 
Lemma \ref{4vectors} implies $$\beta=-1, \alpha=-\bar a\bar c$$
Similarly, the orthogonality of columns $5,6$ gives: $$ \gamma=-ac$$
 Using the formulas for $\alpha,\gamma$ and the orthogonality of columns $3,5$ we obtain:
$$1+\bar a b=0$$
while the orthogonality of columns $3,6$ yields:
$$1-a\bar b=0$$
which is a contradition.

$ $

\noindent\textbf{Case III.}$$H=\HermDiagThreeA$$
The orthogonality of columns 2,4 of $H$ yields:
$$2b+a\bar\alpha+b\bar\beta=0$$
Lemma \ref{4vectors} implies $a\bar\alpha=b\bar\beta=-b$. Thus:
$$\alpha=-a\bar{b}, \beta=-1$$
In particular: $$\alpha\beta=a\bar b$$
 However, using orthogonality of columns 1,4 and Lemma \ref{4vectors} we have: $$\alpha+\beta+\bar b-a=0,\text{ thus: }\{\alpha,\beta\}=\{a,-\bar b\}$$
which implies: $$\alpha\beta=-a\bar b$$
This contradicts $\alpha\beta=a\bar b$.

$ $

\noindent\textbf{Case IV.}$$H=\HermDiagThreeB$$
Orthogonality of columns 4,6 yields:
$$-1-\beta+\alpha\gamma+\beta=0, \text{ thus }\alpha=\bar\gamma$$
Using this together with the orthogonality of columns 1,5 we obtain:
$$2-a-\bar{a}=-2\bar{\gamma}$$
In particular, $\gamma$ has to be real, and since $|a+\bar a|\leq 2$ we must have $\gamma=-1$. Thus, $a+\bar a=0$ so $$a\in\{-i,i\}$$
The sum of the elements of columns 4 and columns 6 is 0. Writing this we obtain: 
$$\bar{\beta}+\beta=0$$
which shows that $\beta\in\{-i,i\}$.
Using now orthogonality of columns 2,3 of $H$ we obtain:
$$b=ac$$
and writing that the sum of the elements of column 4 is 0 we have:
$$ac+c-\beta-1=0$$
By Lemma \ref{4vectors}, we have two possibilities:
$$ac=1, \beta=c\text{ or } c=1,\ \beta=ac$$
Thus, for each choice of $a\in\{-i,i\}$ we have two possible values of $c,\beta$, which uniquely determine the other variables. It is easy to see that all four Butson type matrices we obtain satisfy the hypothesis:

$$ H_1=\begin{pmatrix}
1 &  1           & 1         & 1   &  1     & 1             \cr
1 & -1           & -i        & 1   &  i     & -1      \cr
1 &  i           & -1        & i   &  -i    & -i      \cr
1 &  1           & -i        & -1  & -1     & i    \cr
1 &  -i          & i         & -1  &  1     & -1  \cr
1 &  -1          & i         & -i  & -1     & 1             \cr
\end{pmatrix}$$
$$ H_2=\begin{pmatrix}
1 &  1           & 1         & 1   &  1     & 1             \cr
1 & -1           & i        & 1   &  -i     & -1      \cr
1 &  -i           & -1        & -i   &  i    & i      \cr
1 &  1           & i        & -1  & -1     & -i    \cr
1 &  i          & -i         & -1  &  1     & -1  \cr
1 &  -1          & -i         & i  & -1     & 1             \cr
\end{pmatrix}$$
$$H_3=\begin{pmatrix}
1 &  1           & 1         & 1   &  1     & 1             \cr
1 & -1           & -i        & -i  &  i     & i      \cr
1 &  i           & -1        & 1   & -i     & -1      \cr
1 &  i           & 1         & -1  & -1     & -i   \cr
1 &  -i          & i         & -1  &  1     & -1  \cr
1 &  -i          & -1        & i   & -1     &  1             \cr
\end{pmatrix}$$
$$H_4=\begin{pmatrix}
1 &  1           & 1         & 1   &  1     & 1             \cr
1 & -1           & i        & i  &  -i     & -i      \cr
1 &  -i           & -1        & 1   & i     & -1      \cr
1 &  -i          & 1         & -1  & -1     & i   \cr
1 &  i          & -i         & -1  &  1     & -1  \cr
1 &  i          & -1        & -i   & -1     &  1             \cr
\end{pmatrix}$$
However, we will see that these matrices are in fact equivalent with a certain matrix from the one-parameter family we find in the next case. Thus, it is not necessary to include them in the classification. 

\noindent\textbf{Case V.}$$H=\HermDiagThreeE$$
We show that in this case there exists a one-parameter family of solutions. To obtain the answer in the form given in the statement of the lemma, let us change variables: $a=x, b=-\bar{y}, c=\bar{t},\gamma=z$.  Thus: $$H=\HermDiagThreeEb$$
Since columns 4,5 are orthogonal, we have $yz+\beta \bar z=0$, thus:
$$\beta=-xyz$$
Similarly, the orthogonality of columns 4,6 yields:
$$\alpha=\bar{x}\bar z t$$
Using the formula for $\beta$ in the orthogonality of columns 3,6 we obtain:
$$1+\bar x\bar y-xyz\bar t-\bar t z=0$$
Equivalently: $1+\bar x\bar y=(xy+1)z\bar t$ . Using $1+\bar x\bar y=1+\frac{1}{x}\frac{1}{y}=\frac{1+xy}{xy}$, we obtain:
$$(xy+1)(\frac{1}{xy}-z\bar{t})=0$$
We will assume that $xy\not=\pm 1$. We treat the case $xy=\pm 1$ at the end of the proof. Simplifying by $(xy+1)$ it follows $\frac{1}{xy}-z\bar{t}=0$, thus:
$$z=\bar x \bar y t$$
Using the orthogonality of columns 1,5 we obtain: $$2-x-\bar t+\alpha+\bar z=0$$ and substituting $\alpha,z$ we have: $2-x-\bar t+\bar{x}(xy\bar t) t+xy\bar t=0$. Equivalently, $2-x+y=(1-xy)\bar t$. Since $xy\not=1$ we obtain:
$$\bar t=\frac{2-x+y}{1-xy}$$
which implies:
$$t=\frac{2-\frac{1}{x}+\frac{1}{y}}{1-\frac{1}{xy}}=\frac{2xy+x-y}{xy-1}$$
Since $t\bar t=|t|^2=1$, it follows:
$$\frac{2xy+x-y}{xy-1}\cdot\frac{2-x+y}{1-xy}=1$$
Equivalently: $$(y^2-2y-1)x^2+2(y^2+2y+1)x-(y^2+2y-1)=0$$
Since $y^2-2y-1=0$ does not have solutions of absolute value 1, we must have $y^2-2y-1\not=0$. Solving the above equation for $x$ we obtain:
$$x_1=\frac{-(y^2+2y+1)-\sqrt{2}\sqrt{y^4+2y^3+2y+1}}{y^2-2y-1}$$
$$x_2=\frac{-(y^2+2y+1)+\sqrt{2}\sqrt{y^4+2y^3+2y+1}}{y^2-2y-1}$$
where the square root denotes the principal value of the complex power function $z\rightarrow z^{\frac{1}{2}}$. We need to check if these solutions have absolute value 1 when $|y|=1$. Consider first the case $|x_1|=1$. Denote $\delta=2(y^4+2y^3+2y+1)$. Since $|y|=1$, we have:
$$\bar\delta=2(\frac{1}{y^4}+2\frac{1}{y^3}+2\frac{1}{y}+1)=\frac{1}{y^4}\delta$$
Thus: 
$$1=x_1\bar x_1=\frac{-(y^2+2y+1)-\sqrt{\delta}}{y^2-2y-1}\cdot\frac{-(\bar y^2+2\bar y+1)-\sqrt{\frac{1}{y^4}\delta}}{\bar y^2-2\bar y-1}$$
Depending on $y$, there are two possibilities: $$\sqrt{\frac{1}{y^4}\delta}=-\frac{1}{y^2}\sqrt{\delta}\text{ or }\sqrt{\frac{1}{y^4}\delta}=\frac{1}{y^2}\sqrt{\delta}$$
In the first case, we obtain:
$$\frac{-(y^2+2y+1)-\sqrt{\delta}}{y^2-2y-1}\cdot\frac{-(\bar y^2+2\bar y+1)+\frac{1}{y^2}\sqrt{\delta}}{\bar y^2-2\bar y-1}=1$$
which, after substituting $\bar y=\frac{1}{y}$, becomes:
$$\frac{-(y^2+2y+1)-\sqrt{\delta}}{y^2-2y-1}\cdot\frac{-( y^2+2 y+1)+\sqrt{\delta}}{1-2y-y^2}=1$$
Thus:
$$(-(y^2+2y+1)-\sqrt{\delta})(-(y^2+2y+1)+\sqrt{\delta})=(y^2-2y-1)(1-2y-y^2)$$
Equivalently:
$$(y^2+2y+1)^2-\delta=(y^2-2y-1)(1-2y-y^2)$$
It is immediate to check that this identity holds true for every complex number $y$.

\noindent We now show that the case $\sqrt{\frac{1}{y^4}\delta}=\frac{1}{y^2}\sqrt{\delta}$ leads to a contradiction. We may assume $\delta\not=0$, since we may consider $\delta=0$ as part of the first case. By doing a similar computation, from $|x_1|=1$ we obtain:
$$(-(y^2+2y+1)-\sqrt{\delta})(-(y^2+2y+1)-\sqrt{\delta})=(y^2-2y-1)(1-2y-y^2)$$
However, since we showed that $$(-(y^2+2y+1)-\sqrt{\delta})(-(y^2+2y+1)+\sqrt{\delta})=(y^2-2y-1)(1-2y-y^2)$$ for every $y$, this yields:
$$(-(y^2+2y+1)-\sqrt{\delta})(-(y^2+2y+1)-\sqrt{\delta})=(-(y^2+2y+1)-\sqrt{\delta})(-(y^2+2y+1)+\sqrt{\delta})$$
which implies $$-(y^2+2y+1)-\sqrt{\delta}=-(y^2+2y+1)+\sqrt{\delta}\text{, thus }\delta=0$$
Cancelation of $(-(y^2+2y+1)-\sqrt{\delta})$ was possible, since $x_1\not=0$. We thus obtained a contradiction with the assumption $\delta\not=0$.

\noindent We have shown that: $$|x_1|=1\text{ if and only if }\sqrt{\frac{1}{y^4}\delta}=-\frac{1}{y^2}\sqrt{\delta}$$ We now need to find for what values of $y$ this holds true. 

\noindent Let $y=e^{i\theta}$, where $\theta\in (-\pi,\pi]$ is the principal value of the argument of $y$. We have:
$$\delta=2(y^4+2y^3+2y+1)=2y^2((y^2+y^{-2})+2(y+y^{-1}))=4e^{2i\theta}(cos(2\theta)+2cos(\theta))$$

$ $

\noindent Denote $r=cos(2\theta)+2cos(\theta)$. The equality $\sqrt{\frac{1}{y^4}\delta}=-\frac{1}{y^2}\sqrt{\delta}$ becomes:
$$\sqrt{re^{-2i\theta}}=-e^{-2i\theta}\sqrt{re^{2i\theta}}$$
Denote $\phi=Arg(y^2)\in (-\pi,\pi]$. The previous equality becomes: $$e^{i\phi}\sqrt{re^{-i\phi}}=-\sqrt{re^{i\phi
}}$$
If $r>0$, we have $\sqrt{re^{i\phi}}=\sqrt{r}e^{i\phi/2}$ and $\sqrt{re^{-i\phi}}=\sqrt{r}e^{-i\phi/2}$, thus the equality we want can not hold. If $r<0$, after simplifying by $\sqrt{-r}$ the equality becomes:
$$e^{i(\phi+\pi)}\sqrt{e^{-i(\phi+\pi)}}=\sqrt{e^{i(\phi+\pi)}}$$
which clearly holds true.
 
\noindent We thus need to find the values of $\theta$ such that $$r=cos(2\theta)+2cos(\theta)\leq 0$$ By denoting $p=cos(\theta)$, the inequality becomes $$2p^2+2p-1\leq 0$$ which holds true for $p\in [\frac{-1-\sqrt{3}}{2},\frac{-1+\sqrt{3}}{2}]$. Since $p\in [-1,1]$, we obtain $cos(\theta)\in [-1,\frac{-1+\sqrt{3}}{2}]$, hence $$\theta\in [-\pi,-arcos(\frac{-1+\sqrt{3}}{2})]\cup [arcos(\frac{-1+\sqrt{3}}{2}),\pi]$$
We have thus obtained a one-parameter family of matrices, which can be easily checked to be Hadamard:
$$H(\theta)=\begin{pmatrix}
1 & 1 & 1 & 1 & 1 & 1 \cr
1 & -1 & \frac{1}{x} &-y  &-\frac{1}{x}  &y   \cr
1 & x & -1 & t & -t  & -x  \cr
1 &-\frac{1}{y}  & \frac{1}{t} & -1  & \frac{1}{y} &-\frac{1}{t}   \cr
1 &-x  & -\frac{1}{t} & y & 1 &\frac{1}{z}  \cr
1 & \frac{1}{y} &-\frac{1}{x}  &-t  & z & 1 \cr
\end{pmatrix}$$
where:
$$y=exp(i\theta),\text{ } z=\frac{1+2y-y^2}{y(-1+2y+y^2)}$$
$$x=\frac{1+2y+y^2-\sqrt{2}\sqrt{1+2y+2y^3+y^4}}{1+2y-y^2}$$
$$t=\frac{1+2y+y^2-\sqrt{2}\sqrt{1+2y+2y^3+y^4}}{-1+2y+y^2}$$
A similar analysis for $|x_2|=1$ leads to another one-parameter family of solutions, for the same interval of values of $\theta$:
$$H'(\theta)=\begin{pmatrix}
1 & 1 & 1 & 1 & 1 & 1 \cr
1 & -1 & \frac{1}{x} &-y  &-\frac{1}{x}  &y   \cr
1 & x & -1 & t & -t  & -x  \cr
1 &-\frac{1}{y}  & \frac{1}{t} & -1  & \frac{1}{y} &-\frac{1}{t}   \cr
1 &-x  & -\frac{1}{t} & y & 1 &\frac{1}{z}  \cr
1 & \frac{1}{y} &-\frac{1}{x}  &-t  & z & 1 \cr
\end{pmatrix}$$
where:
$$y=exp(i\theta),\text{ } z=\frac{1+2y-y^2}{y(-1+2y+y^2)}$$
$$x=\frac{1+2y+y^2+\sqrt{2}\sqrt{1+2y+2y^3+y^4}}{1+2y-y^2}$$
$$t=\frac{1+2y+y^2+\sqrt{2}\sqrt{1+2y+2y^3+y^4}}{-1+2y+y^2}$$
However, it is easy to check that $H'(\theta)$ is equivalent to $H(\theta)$:
$$P_1D_1H(\theta)D_2P_2=H'(\theta)$$
where $D_1$ is the unitary diagonal matrix:
$$\begin{pmatrix} 1 & 0 & 0 & 0 & 0 & 0 \cr 0 & -1 & 0 & 0 & 0 & 0 \cr 
0 & 0 & - \frac{-1 - 2\,y + y^2}
     {1 + 2\,y + y^2 - {\sqrt{2}}\,{\sqrt{1 + 2\,y + 2\,y^3 + y^4}}} 
  & 0 & 0 & 0 \cr 0 & 0 & 0 &
    -y & 0 & 0 \cr 0 & 0 & 0 & 0 & \frac{-1 - 2\,y + y^2}
   {1 + 2\,y + y^2 - {\sqrt{2}}\,{\sqrt{1 + 2\,y + 2\,y^3 + y^4}}} & 0 
\cr 0 & 0 & 0 & 0 & 0 & y \cr 
\end{pmatrix}$$
and $D_2$ is the unitary diagonal matrix:
$$\begin{pmatrix} -1 & 0 & 0 & 0 & 0 & 0 \cr 0 & 1 & 0 & 0 & 0 & 0 \cr 
0 & 0 & \frac{1 + 2\,y + y^2 -
     {\sqrt{2}}\,{\sqrt{1 + 2\,y + 2\,y^3 + y^4}}}{-1 - 2\,y + y^2} & 0 
& 0 & 0 \cr 0 & 0 & 0 & \bar y & 0 & 0 \cr 0 & 0 & 0 & 0 & \frac{-1 - 2\,y - y^2 + 
{\sqrt{2}}\,{\sqrt{1 + 2\,y + 2\,y^3 + y^4}}}
   {-1 - 2\,y + y^2} & 0 \cr 0 & 0 & 0 & 0 & 0 & -\bar y  \cr 
\end{pmatrix}$$
and $P_1,P_2$ are permutation matrices:
$$P_1=\begin{pmatrix}
   0 & 1 & 0 & 0 & 0 & 0 \cr 0 & 0 & 0 & 0 & 0 & 1 \cr 0 & 0 & 0 & 0 & 1 
& 0 \cr 1 & 0 & 0 & 0 & 0 & 0 \cr 0 &
   0 & 0 & 1 & 0 & 0 \cr 0 & 0 & 1 & 0 & 0 & 0 \cr 
\end{pmatrix}, P_2=\begin{pmatrix}
   0 & 0 & 0 & 1 & 0 & 0 \cr 1 & 0 & 0 & 0 & 0 & 0 \cr 0 & 0 & 0 & 0 & 0 
& 1 \cr 0 & 0 & 0 & 0 & 1 & 0 \cr 0 &
   0 & 1 & 0 & 0 & 0 \cr 0 & 1 & 0 & 0 & 0 & 0 \cr 
\end{pmatrix}$$
We should also mention here that the four matrices $H_1,H_2,H_3,H_4$ found in case 4, which are easy to check to be equivalent, are in fact equivalent to the matrix $H(\pi/2)$ which is part of the one-parameter family. Indeed, we have:
$$H(\pi/2)=P_1DH_3P_2$$
where:
$$D=\begin{pmatrix} 1 & 0 & 0 & 0 & 0 & 0\\
0 & i & 0 & 0 & 0 & 0\\
0 & 0 & -1 & 0 & 0 & 0\\
0 & 0 & 0 & 1 & 0 & 0\\
0 & 0 & 0 & 0 & -i & 0\\
0 & 0 & 0 & 0 & 0 & -1\\
\end{pmatrix}$$
$$P_1=\begin{pmatrix} 1 & 0 & 0 & 0 & 0 & 0 \\
0& 1 & 0 & 0 & 0 & 0\\
0 & 0 & 1 & 0 & 0 & 0\\
0 & 0 & 0 & 0 & 1 & 0 \\
0 & 0 & 0 & 0 & 0 & 1 \\
0 & 0 & 0 & 1 & 0 & 0 \\
\end{pmatrix}, P_2=\begin{pmatrix} 0 & 0 & 0 & 0 & 0 & 1\\
0 & 0 & 0 & 1 & 0 & 0\\
1 & 0 & 0 & 0 & 0 & 0 \\
0 & 0 & 1 & 0 & 0 & 0\\
0 & 0 & 0 & 0 & 1 & 0\\
0 & 1 & 0 & 0 & 0 & 0\\
\end{pmatrix}$$
To end the proof, we still need to show that $xy=\pm 1$ leads to no solutions. 
We first consider the case $xy=-1$, thus $y=-\bar x$ and $\beta=-xyz=z$. Since the sum of rows 4,5 is 0, we have:
$$2-1+\bar\alpha+\alpha+1+z+z=0$$
thus:
$$1+Re(\alpha)=-z$$
which implies that $z$ is real, and since $|z|=1$ and $|Re(\alpha)|\leq 1$ we must have $z=1$, $Re(\alpha)=0$.
However, since the sum of the elements of row 5 is 0:
$$3=\bar x+t-\bar \alpha$$
and the triangle inequality implies $\alpha=-1$, contradicting $Re(\alpha)=0$.

\noindent Consider now the case $xy=1$, so $y=\bar x$ and $\beta=-xyz=-z$. Summing up columns 5,6 we obtain:
$$4=x+t-\bar\alpha-z$$
and the triangle inequality shows $x=t=1$, $\alpha=z=-1$. However this implies that the sum of elements of column 5 is $1-1-1-1+1-1=-2\not=0$, contradiction.

\end{proof}
\begin{lem}
Let $H\in M_6(\mathbb C)$ be a self-adjoint, dephased, complex Hadamard matrix. Then the diagonal of $H$ can not be $(1,-1,-1,-1,-1,1)$. 
\end{lem}
\begin{proof}
Reasoning as in the previous lemmas, we may assume:  
$$H=\HermDiagFour$$
Using the orthogonality of columns 3,5 we have:
$$-2y+x\alpha+z\bar\gamma=0$$
Lemma \ref{4vectors} implies $$ y=\alpha x, z= \gamma \alpha x$$
Using this and the orthogonality of columns 4,6 we obtain:
$$\alpha \gamma=0$$
which is not possible since $|\alpha|=|\gamma|=1$.

\end{proof}

\begin{lem}
Let $H\in M_6(\mathbb C)$ be a self-adjoint, dephased, complex Hadamard matrix. Then the diagonal of $H$ can not be $(1,-1,-1,-1,-1,-1)$.
\end{lem}
\begin{proof}
We may assume
$$H=\HermDiagFive$$
Using the fact that columns 3,5 are orthogonal, as in the previous lemma, we obtain:
$$ y=\alpha x, z= \gamma \alpha x$$
Now the orthogonality of columns 4,6 together with lemma \ref{4vectors} yields:
$$\beta=\alpha\gamma=\bar{x}z$$
Using the expression for $\beta$ and the fact that columns $3,6$ are orthogonal we obtain:
$$1-\bar a b-z+x\beta+y\gamma-z=0$$
Since $x\beta=z$ and $y\gamma=z$, it follows
$$1-\bar a b=0$$
Thus $b=a$. But in this case the inner product of columns 3,4 is 
$$2-2x+y\bar \alpha+z\bar\beta=2-2x+2x=2\not=0$$
which contradicts the fact that they are orthogonal.\end{proof}

\noindent This ends the proof of Theorem \ref{main}. 
\begin{obs} The Butson type matrices $H_1,H_2,H_3,H_4$, which are equivalent to $H(\pi/2)$, are also equivalent to the matrix $D_6$ from the catalogue \cite{TZ}. This shows that besides the affine family through $D_6$, exhibited in \cite{TZ}, there also exists a one-parameter non-affine family containing $D_6$.
\end{obs}

\section{More Complex Hadamard Matrices}
$ $

We present in this section some new examples of complex Hadamard matrices of small orders. These examples were found using Mathematica, by searching for the local minimum of a function $f(H)$ encoding the conditions a matrix $H$ needs to satisfy to be Hadamard, i.e. $f(H)=0$ if and only if $H$ is Hadamard. 

Trying to find such general solutions leads to numerous numerical results that are difficult to interpret formally. However, it turns out that asking for extra symmetries for $H$, such as $H$ hermitian, often yields some clear algebraic results. 

Inspired by our work on the classification of $6\times 6$ Hadamard matrices, we tried to fix certain diagonal entries of $H$, making them $1$ or $-1$. Such conditions seem to be very strong, leading to solutions that are easy to interpret algebraically.

$ $

\noindent\textbf{New examples for n=9.}

$ $

 $H_{9}=$\CircNine, $y=\frac{1- i\sqrt{15}}{4}$

It can be easily checked that $H_{9}$ satisfies the "span condition" of \cite{Ni}, thus it is isolated among Hadamard matrices. In particular, it is not equivalent to any of the matrices from the 4-parameter family $F_9^{(4)}$ described in \cite{TZ}. Also, $H_9$ is not equivalent to any of the known circulant solutions, since its entries belong to $\mathbb Q[\frac{1- i\sqrt{15}}{4}]$.

The next matrix was obtained by searching for the local minimum with a fixed diagonal:

$BN_{9}=$ \nineten, $\epsilon=e^{2\pi i/10}$

The matrix $BN_{9}$ has defect 2, in the sense of \cite{TZ}, thus it might be part of a family of non-equivalent Hadamard matrices.  

$ $

\noindent\textbf{New examples for n=10.}   

$ $

$BN_{10}=$\liketau , $a=\frac{-1+i\sqrt{15}}{4}$

This matrix was found by numerical search of the local minimum, with the constraint that all diagonal entries be 1. It satisfies the span condition and thus it is isolated among complex Hadamard matrices. 
 
$ $

\noindent\textbf{New examples for n=11.} 

$ $ 

$N_{11},N'_{11}=$\romeo

\noindent where $x= \frac{3}{4}-i \frac{\sqrt{7}}{4}$ and
$\begin{pmatrix} a \\ b \\ y \\ \end{pmatrix} = \begin{pmatrix} -x \\ x\\ -x^2  \end{pmatrix}$ or 
$\begin{pmatrix} a \\ b \\ y \\ \end{pmatrix} = \begin{pmatrix} 1 \\ -1 \\ \bar{x} \end{pmatrix}$.  
Both these matrices are isolated.

\bibliographystyle{amsalpha}

\begin{thebibliography}{10}
\bibitem[Bj]{Bj}G. Bjorck, \emph{Functions of modulus 1 on $\mathbb{Z}_n$, whose Fourier transforms have constant modulus, and "cyclic n-roots"}, Recent Advances in Fourier Analysis and its applications, NATO Adv. Sci. Inst. Ser. C: Math. Phys. Sci., Kluwer \textbf{315} (1990), 131-140.


\bibitem[CH]{CH} C. H. Cooke, I. Heng, \emph{Error correcting codes associated with complex Hadamard matrices}, Appl. Math. Lett. \textbf{11} (1998), 77

\bibitem[Di]{Di}P. Dita, \emph{Some results on the parametrization of complex Hadamard matrices}, J. Phys. A, \textbf{37} (2004) no. 20, 5355-5374

\bibitem[GHJ]{GHJ}F. M. Goodman, P.de la Harpe and V.F.R.Jones,\emph{Coxeter Graphs and Towers of Algebras}, Math. Sciences Res. Inst. Publ. Springer Verlag 1989

\bibitem[H]{Haagerup}U. Haagerup, \emph{Orthogonal maximal abelian *-subalgebras of the $n\times n$ matrices and cyclic $n$-roots}, Operator Algebras and Quantum Field Theory (ed. S.Doplicher et al.), International Press (1997),296-322

\bibitem[HJ]{HJ}P. de la Harpe and V. F. R. Jones, \emph{Paires de sous-algebres semi-simples et graphes fortement reguliers}, C.A. Acad. Sci. Paris \textbf{311}, serie I (1990), 147-150

\bibitem[Jo1]{Jo1}V. F. R. Jones, \emph{Index for subfactors}, Invent. Math \textbf{72} (1983),
  1--25.

\bibitem[Jo2]{Jo2}V. F. R. Jones, \emph{Planar Algebras I}, math.QA/9909027

\bibitem[MRS]{MRS}M. Matolcsi, J. Reffy, F. Szollosi, \emph{Constructions of complex Hadamard matrices via tiling abelian groups}, preprint quant-ph/0607073

\bibitem[MW]{MW}A. Munemasa and Y. Watatani, \emph{Orthogonal pairs of *-subalgebras and Association Schemes} C.R. Acad. Sci. Paris \textbf{314}, serie I (1992), 329-331.

\bibitem[Ni]{Ni}R. Nicoara, \emph{A finiteness result for commuting squares of matrix algebras}, Journal of Operator Theory 2006, math.OA/0404301 

\bibitem[Pe]{Petrescu}M. Petrescu, \emph{Existence of continuous families of complex Hadamard matrices of certain prime dimensions and related results}, PhD thesis, University of California Los Angeles, 1997.

\bibitem[Po1]{Po1}S. Popa, \emph{Classification of subfactors : the reduction to commuting squares}, Invent. Math., \textbf{101}(1990),19-43


\bibitem[Po2]{Po2}S. Popa, \emph{Othogonal pairs of *-subalgebras in finite von Neumann algebras}, J. Operator Theory \textbf{9}, 253-268 (1983)

\bibitem[TZ]{TZ}
W. Tadej and K. Zyczkowski,\emph{ A concise guide to complex Hadamard matrices}
Open Systems \& Infor. Dyn.,\textbf{13}(2006), 133-177, quant-ph/0512154 

\bibitem[T]{T} T. Tao, \emph{Fuglede's conjecture is false in 5 and higher dimensions}, Math. Res. Letters \textbf{11} (2004), 251

\bibitem[We]{We}R. F. Werner, \emph{All teleportation and dense coding schemes}, J.Phys.A, \textbf{34} (2001), 7081-7094
\end{thebibliography}

\end{document}